\documentclass{amsart}
\usepackage{amsmath,amssymb,amsthm}
\theoremstyle{definition}

\theoremstyle{remark}

\numberwithin{equation}{section}

\begin{document}

\title{A counterexample on spectra of zero patterns}

\author{Yaroslav Shitov}
\address{National Research University Higher School of Economics, 20 Myasnitskaya Ulitsa, Moscow 101000, Russia}
\email{yaroslav-shitov@yandex.ru}


\begin{abstract}
An $n\times n$ \textit{zero pattern} $S$, which is a matrix with entries $*$ and $0$, is called \textit{spectrally arbitrary} with respect to a field $\mathbb{F}$ if any monic polynomial $f$ of degree $n$ can be realized as the characteristic polynomial of a matrix obtained from $S$ by replacing the $*$'s with non-zero elements of $\mathbb{F}$. We construct an $n\times n$ zero pattern that is spectrally arbitrary with respect to $\mathbb{C}$ and has $2n-1$ nonzero entries. 
\end{abstract}

\maketitle

One of the intriguing questions in the theory of incomplete matrices is as follows. What is the smallest number $k$ such that there exists an $n\times n$ zero pattern with $k$ nonzero entries which is spectrally arbitrary over $\mathbb{R}$? It was proved in~\cite{DJOD} that $k$ has to be at least $2n-1$ in the above question, and we note that the proof works even if $\mathbb{R}$ gets replaced by any other infinite field. On the other hand, many spectrally arbitrary $n\times n$ zero patterns with $2n$ non-zero entries are known (see~\cite{GS}), so the optimal value of $k$ is at most $2n$. The statement that this optimal value is in fact $2n$ has become known\footnote{Strictly speaking, the $2n$ conjecture is the same question but asked for the sign patterns instead of zero patterns.} as the $2n$ \textit{conjecture} (see~\cite{alot}), and it attracts a significant amount of attention in the contemporary linear algebra community (\cite{BMOD, COD, DJOD, McY}).

This paper is a piece of evidence against this conjecture. Although we were not able to disprove it so far (and we explain why in the end of the paper), we present a counterexample to its complex analogue. Also, our result allows us to answer several questions asked by McDonald and Yielding in~\cite{McY}.

Now we proceed with a counterexample. Throughout our paper, we denote by $x_1,\ldots,x_8$ a family of variables that are allowed to take any complex value except zero. We consider the matrix $X(x_1,x_2,x_3,x_4,x_5,x_6,x_7,x_8)$ defined as
$$\begin{pmatrix}
x_1&1&0&0&0&0&0&0\\
x_7&x_2&1&0&0&0&0&0\\
0&x_3&0&1&0&0&0&0\\
0&0&0&0&1&0&0&0\\
0&0&0&0&0&1&0&0\\
0&x_8&0&0&0&0&1&0\\
0&0&0&0&0&0&0&1\\
x_6&0&0&x_5&0&x_4&0&0
\end{pmatrix},$$
and we denote its zero pattern by $S$. It is easy to check that the matrix $X(1,-1,1,1,-1,1,-2,1)$ is nilpotent, so we can realize $t^8$ as the characteristic polynomial of a matrix with zero pattern $S$. The polynomial $(x-1)^8$ is a bit harder to realize: One needs to take $x_1=1737/848$, $x_2=5047/848$, $x_3=-4452/193$, $x_4=35/4$, $x_5=2/7$, $x_6=25/2$, $x_7=1007374319/138787072$, $x_8=-1325/7$, turn on the computer, and check that the matrix $X$ defined this way has a desired characteristic polynomial.

How to describe all the characteristic polynomials realized by $S$? First of all, we note that any matrix with pattern $S$ can be reduced to the form $X$ by conjugating it with a diagonal matrix. Since a pair of similar matrices have the same characteristic polynomial, we can restrict our attention to the matrices of the form $X$. We use the computer again and compute $\varphi=\det (tI-X)$, which is a monic polynomial $t^8+\varphi_7 t^7+\ldots+\varphi_0$ with coefficients $\varphi_i$ in $\mathbb{C}[x_1,\ldots,x_8]$. As we look into the result of the computation, we note that $\varphi_4$ is a multiple of $\varphi_7$. In particular, $\varphi_7$ cannot be zero unless $\varphi_4$ is zero, so $S$ is not a spectrally arbitrary pattern.

Despite this fact, a lot of polynomials can be realized by the pattern $S$. To see this, we consider new variables $\tau_0,\ldots,\tau_7$ and find the simultaneous solution of the equations $\varphi_i=\tau_i$ for $x_1,\ldots,x_8$. This may sound as a hard task, but the computer comes up with the solution immediately, --- this becomes possible thanks to a carefully selected pattern $S$. The resulting values of the $x_i$'s are all rational functions in the $\tau_i$'s, so the polynomial $t^8+\tau_7t^7+\ldots+\tau_0$ can be realized by $S$ unless one of these rational functions has a vanishing numerator or denominator. Taking the LCM of all these numerators and denominators, we get a polynomial $\pi(\tau_0,\ldots,\tau_7)$ which can vanish only if $t^8+\tau_7t^7+\ldots+\tau_0$ cannot be realized by $S$.

Now we have learned everything we need about the possible characteristic polynomials of $X$, and we turn our attention to its \textit{spectra}. For any family $\sigma$ of eight complex numbers, we define $s_i$ to be $(-1)^i$ times the $i$th elementary symmetric polynomial of $\sigma$ and define $\psi(\sigma)=\pi(s_8,s_7,\ldots,s_1)$. One can determine the total degree\footnote{Of course, we do not need to calculate $\psi$ to do this. Instead, we compute the maximal weight of the monomials of $\varphi$ with respect to the function $(\tau_0,\ldots,\tau_7)\to(8,7,\ldots,1)$.} of $\psi$ (which equals $94$), and Vieta's formulae show that $\sigma$ is the spectrum of a matrix of the form $X$ whenever $\psi$ does not vanish.

Since $\psi$ is a polynomial in eight variables and has total degree $94$, every set of $94+8=102$ distinct complex numbers has a subset $\sigma$ of eight elements which satisfies $\psi(\sigma)\neq0$. Therefore, if a family $U$ of $708$ complex numbers contains at least $102$ distinct elements, then it has a subset realizable as the spectrum of a matrix of the form $X$. Otherwise, we use the pigeonhole principle and conclude that some number $c$ repeats in $U$ at least eight times. As explained above, $(x-c)^8$ is realizable as the characteristic polynomial of $X$, so $U$ does anyway contain a subfamily $V$ realizable as the spectrum of a matrix $M$ with pattern $S$.

Now we define $D_{2m}$ to be the block-diagonal matrix consisting of the $m$ blocks equal to the $2\times2$ matrices of all $*$'s. It is easy to see that $D_{2m}$ is spectrally arbitrary, so we can find a matrix $M'$ with pattern $D_{700}$ and spectrum $U\setminus V$. Now we see that the matrix $\operatorname{diag}(M,M')$ has spectrum $U$ and pattern $\operatorname{diag}(S,D_{700})$. Therefore, $\operatorname{diag}(S,D_{700})$ is a $708\times 708$ zero pattern which has $1415$ nonzero elements and is spectrally arbitrary with respect to $\mathbb{C}$.

As said above, our result answers two questions asked in~\cite{McY}. First, we have constructed an $n\times n$ zero pattern which has $2n-1$ nonzero entries and is spectrally arbitrary with respect to $\mathbb{C}$. Secondly, we get an example of zero patterns $A,B$ such that $\operatorname{diag}(A, B)$ is spectrally arbitrary but $A$ is not. We note in passing that an argument similar to the proof of item (3) of Theorem~11 in~\cite{mysign} would allow us to get a refined version of this result. Namely, we would be able to show that $\operatorname{diag}(P, P)$ can be a spectrally arbitrary pattern even if $P$ is not spectrally arbitrary.

We conclude our paper with several thoughts on the real version of the $2n$ conjecture. My attempts to disprove it were not successful, and the main obstacle was the fact that the set of real irreducible polynomials is much richer than its complex counterpart. A brief examination of our counterexample allows one to prove the following sufficient condition for $S$ to be a diagonal block of a block-diagonal pattern spectrally arbitrary over $\mathbb{C}$. Namely, this happens if the set of all characteristic polynomials allowed by $S$ contains $t^n$, $(t-1)^n$, and a generic monic polynomial of degree $n$. Unfortunately, this condition is not sufficient for real patterns, and the reason lies in the existence of degree-two irreducible polynomials over $\mathbb{R}$. The set of such polynomials remains positive-dimensional even if the spectra are considered up to scaling, which makes it hard to produce families of matrices that depend on $n$ parameters only and allow all such polynomials. For instance, there exist (finitely many) values of $a\in(-2,2)$ for which $(t^2+at+1)^4$ cannot be the spectrum of the matrix $X$ as above. In particular, this happens if
$$a=0\,\,\,\mbox{or}\,\,\,a=\sqrt{\frac{\sqrt{15}-3}{3}}$$
and makes it impossible for $S$ to be a diagonal block of a block-diagonal pattern spectrally arbitrary over $\mathbb{R}$. However, we believe that a counterexample for the real $2n$ conjecture can be found with a more extensive search.

\end{document}